\documentclass[11pt]{article}
\usepackage[english]{babel}
\usepackage{booktabs}
\usepackage{float}

\usepackage{caption}

\usepackage{amsmath}
\usepackage{amsfonts}
\usepackage{mathtools}
\usepackage{graphicx}
\usepackage[colorlinks=true, allcolors=blue]{hyperref}
\usepackage[normalem]{ulem}
\usepackage{soul}
\usepackage{xcolor}
\usepackage{amsthm}
\usepackage{blindtext}

\usepackage{tikz-cd}
\newsavebox{\pullback}
\sbox\pullback{%
\begin{tikzpicture}%
\draw (0,0) -- (1ex,0ex);%
\draw (1ex,0ex) -- (1ex,1ex);%
\end{tikzpicture}}

\title{Fox-Milnor condition for concordant knots in homology 3-spheres}
\author{Bao Vuong }
\date{April 2025}

\begin{document}

\maketitle

This paper will show that the Alexander polynomial of a knot, which is of slice type in an oriented homology 3-sphere, obeys the Fox-Milnor polynomial condition. A relation between Alexander polynomial of concordant knots in an oriented homology 3-sphere is established.

\medskip
\noindent\textbf{Note added in this version.}
In the published version, the divisibility in Theorem A  and Theorem B was stated in the polynomial ring. The correct statement is in the field of fractions; equivalently, a factor \(q(t)q(t^{-1})\) on the left-hand side must be included. With this correction the proof is unchanged. I thank Mark Powell for pointing this out and for bringing to my attention Chapter 15 of \emph{The foundations of four-manifold theory in the topological category \cite{MarcPowell} }, available at
\url{https://nyjm.albany.edu/m/}, where a related Fox--Milnor theorem for knots in integral homology 3-spheres bounding integral homology 4-balls is proved. That work appeared prior to this paper and should have been cited.

%Note added in v2 (3 Sept 2026): A minor correction has been made to the proof of Theorem A and Theorem B. The theorems and all results mostly remain unchanged. Precisely, for the equation, matching the two Alexander polynomials, must come along with an element of the field of fractions of  $\mathbb{Z}[t^{\pm}]$. A reference to the prior work of Stefan Friedl, Matthias Nagel, Patrick Orson, and Mark Powell \cite{MarcPowell} has been added, following post-publication correspondence.

\section{Introduction}
The classical knot theory considers the problem of placement of a circle into three dimensional sphere, that is to classify topological embedding up to ambient isotopy. Another essential problem is classification of knots up to cobordism. The latter problem leads to research and building cobordism invariant. One of such invariant was explored by Fox and Milnor, announced in 1957 and published in their classical paper \cite{FoxMilnor}. In the paper,  Fox and Milnor proved a condition for the Alexander polynomial for a knot being a slice type. In this paper, we consider knots in an oriented homology 3-sphere. For a knot in an oriented homology 3-sphere, there are Alexander invariants, associated with the knot in a manner of classical knots (see \cite{Saveliev}, p.~91--93). That is, the geometric essence, generating the Alexander invariants of knots in an oriented homology 3-sphere, is established via the infinite cyclic covering of the knot complement. Particularly for each knot in an oriented homology 3-sphere there is a well-defined Alexander polynomial, related to the knot. The polynomial is an invariant under ambient isotopy of knots in an oriented homology 3-sphere. A concise text about Alexander polynomials of knots in an oriented homology 3-sphere can be found in the paper \cite{AustinRolfsen}. For a special case of Poincare homology sphere $P$, we give an algorithm to compute the Alexander polynomial for a knot or a link in $P$ in a recent work \cite{BaoEvteev}. It is natural to ask if the Fox-Milnor condition on Alexander polynomial holds for the case of concordant knots in an oriented homology 3-sphere.
We will give an affirmative answer by proving the following theorems

{\bf Theorem A}. Let $k_0, k_1$ be concordant knots in an oriented homology 3-sphere $M$. Then the Alexander polynomials of the knots are related by the following equation \footnote{The notation $\Delta_1 \dot{=} \Delta_2$ means $\Delta_1 = \pm t^n \Delta_2$ for some integer $n$}

\[\Delta_{k_0}(t) ~\dot{=} ~\lambda(t)\lambda(1/t) \Delta_{k_1}(t)\]

where $\Delta_{k_0}(t), \Delta_{k_1}(t)$ are the Alexander polynomials in $t$ of the knots $k_0,k_1$ respectively and $\lambda(t)$\footnote{In the published version, this statement was formulated in $\mathbb{Z}[t^{\pm}]$. The correct statement is in the field of fractions of $\mathbb{Z}[t^{\pm}]$, which is equivalent to including the factor \(q(t)q(t^{-1})\) to the left-hand side. See the Note added in this version.} is an element of the field of fractions of  $\mathbb{Z}[t^{\pm}]$.

{\bf Theorem B}. Let $M, M'$ be homological spheres. Let $\mathcal{W}$ be a cobordism between $M$ and $M'$, and the boundary of $\mathcal{W}$ is disjoint union $\partial \mathcal{W} = M \cup M'$. More over the inclusions $M \hookrightarrow \mathcal{W}$ and $M' \hookrightarrow \mathcal{W}$ induce isomorphisms on homology. Let $k$ and $k'$ be knots in $M$ and $M'$ correspondingly. If there exist a concordance $g: S^1 \times I \rightarrow \mathcal{W}$ between $k$ and $k'$. Then the Alexander polynomials of the knots $k$ and $k'$ are related by the following equation

\[\Delta_{k}(t) ~\dot{=} ~\lambda(t)\lambda(1/t) \Delta_{k'}(t)\]

where $\Delta_{k}(t), \Delta_{k'}(t)$ are the Alexander polynomials in $t$ of the knots $k,k'$ respectively and $\lambda(t)$\footnote{See Note added in this version. } is an element of the field of fractions of  $\mathbb{Z}[t^{\pm}]$.

To prove the Theorem A we follow a similar strategy, that of Fox and Milnor in \cite{FoxMilnor} which is in the form of Kreinbihl \cite{Kreinbihl} and Turaev \cite[p.~139]{Turaev}. Before that we need some information about knot complement in the homology 3-sphere $M$ and corresponding concordance complement. Theorem B is a direct generalisation from Theorem A.

\section{Knots in an oriented homology 3-sphere}

A knot $k$ in an oriented homology 3-sphere $M$ is a smooth embedding $k: S^1 \rightarrow M$. Two knots $k_1,k_2$ in the homology 3-sphere $M$ are equivalent if there is an ambient isotopy taking one to the other. A knot is call trivial or an unknot if it bounds a smoothly embedded disk in $M$. We carry a proof of the following well-known fact \footnote{We do not know any text on proof of this elsewhere for except the fact is stated in the book \cite{Saveliev} by Saveliev (in page 91 and exercise 2 in page 42)}

{\bf Lemma 1}. Let $k$ be a knot in an oriented homology 3-sphere $M$. Then the knot complement $M-k$ in $M$ has the homology of a circle, that is $H_*(M-k)=H_*(S^1)$.

{\it Proof.}

We will proceed by using Mayer-Vietoris sequence.
Let $\nu(k)$ be a tubular neighbourhood of the knot $k$ in the homology 3-sphere $M$. The homology 3-sphere $M$ is the union $(M-k) \cup \nu(k)$ of the knot exterior $M-k$ and the tubular neighbourhood $\nu(k)$. The intersection $\nu(k) \cap (M-k)$ is retracted to a torus $T^2$. The tubular neighbourhood $\nu(k)$ is retracted to its core circle $S^1$ and its boundary $\partial \nu(k)$ is homeomorphic to a torus. Note that, the space $M-k$ deform retracts onto closure $X_k$ of $M-\nu(k)$ then $H_{*}(M-k) \cong H_*(M-\nu(k)$.
Since $M-\nu(k)$ is a connected, 3-manifold with boundary for any cellular structure on it, we can collapse $M-k$ down to a 2-dimensional subcomplex by pushing all 3-cells to the boundary. This collapsing do not affect the homology. Hence we get the result of the theorem for all $* \neq 1,2$. Consider the following part of the Mayer-Vietoris sequence:
\begin{align*} 
H_3(T^2) & \rightarrow H_3(\nu(k)) \oplus H_3(M-k) \rightarrow H_3(M) \rightarrow \\ 
\rightarrow H_2(T^2) & \rightarrow H_2(\nu(k)) \oplus H_2(M-k) \rightarrow H_2(M) \rightarrow\\
\rightarrow H_1(T^2) & \rightarrow H_1(\nu(k)) \oplus H_1(M-k) \rightarrow H_1(M)
\end{align*}

that is equivalent to
\begin{align*} 
0 &\rightarrow 0 \oplus 0 \rightarrow \mathbb{Z} \rightarrow \\ 
\rightarrow \mathbb{Z} & \rightarrow 0 \oplus H_2(M-k) \rightarrow 0 \rightarrow\\
\rightarrow \mathbb{Z}^2 & \rightarrow \mathbb{Z} \oplus H_1(M-k) \rightarrow 0 
\end{align*}

The map $\mathbb{Z}^2 \rightarrow \mathbb{Z} \oplus H_1(M-k)$ is an isomorphism by exactness, so we get $H_1(M-k)=\mathbb{Z}$ \footnote{This can be done immediately by Alexander duality (\cite{AustinRolfsen}, page 258)}

Consider the portion containing $H_2(M-k)$

$$0 \rightarrow \mathbb{Z} \xrightarrow{\text{injective}} \mathbb{Z} \xrightarrow{\text{surjective}} H_2(M-k) \rightarrow 0$$

From this short exact sequence we can conclude that $H_2(M-k)$ is finite cyclic.
By Poincar\'e-Lefschetz duality and the Universal Coefficients Theorem we get
$$H_2(M-k) \cong H_2(X_k) \cong H^1(X_k,\partial X_k) \cong \text{Hom}(H_1(X_k,\partial X_k), \mathbb{Z}).$$
As $\text{Hom}(H_1(X_k,\partial X_k), \mathbb{Z})=\mathbb{Z}^{\beta_1(X_k,\partial X_k)}$, where $\beta_1$ is the first Betti number, hence $H_2(M-k)$ is a free abelian group. That means $H_2(M-k)=0$. Q.E.D.

\section{Concordance of knots in an oriented homology 3-sphere}

In classical knot theory, two tame knots (smooth knots, p.l. knots) $k_0$ and $k_1$ in 3-sphere are said to be concordant\footnote{For a survey of classical knot concordance see Livingston \cite{Livingston}} if the knots are the boundary of a smooth embedding of the cylinder $S^1 \times [0,1]$ into $S^3 \times [0,1]$, that is $k_0 \subset S^3 \times \{0\}, k_1 \subset S^3 \times \{1\}$. A knot is called a slice knot if it is concordant to the unknot (compare with definition by Fox-Milnor \cite[p.~258]{FoxMilnor}). We can give an analogous definition for concordant knots in an oriented homology sphere

{\bf Definition.} Let $k_0, k_1$ be knots in an oriented homology 3-sphere $M$. The two knots are smoothly concordant if there is a smooth embedding $f$ of the cylinder $S^1 \times [0,1]$ into $M \times [0,1]$ such that
$$f(S^1 \times \{0\}=k_0 \times \{0\} \subset M \times \{0\}$$
and
$$f(S^1 \times \{1\}=k_1 \times \{1\} \subset M \times \{1\}$$

That is $k_0$ and $k_1$ cobound a smooth embedded cylinder in $M\times I$, where $I=[0,1]$.

A knot $k$ in the homology 3-sphere $M$ is called slice if $k$ concordant to the unknot\footnote{A definition of topological sliceness can be given by the notion of locally-flat embedded disk in a contractible topological 4-manifold (see \cite[Corollary 4]{AustinRolfsen})}.

We compute the homology of concordance complement and determine its boundary. Let $C$  be the embedded cylinder, the image of $S^1 \times I$ under concordance $f$. Let $N(C)$ be a tubular neighbourhood of the cylinder $C$ in $M \times I$.

{\bf Lemma 2}. Let $k_0,k_1$ be concordant knots in an oriented homology 3-sphere $M$. Let $C$ and $N(C)$ be defined as above. The knots $k_0,k_1$ cobound the cylinder $C$. Then concordance complement $(M\times I - N(C))$ has the homology of the circle $S^1$, that is
$H_*(M\times I-N(C))= H_*(S^1)$.

{\it Proof.} Lemma 2 is a special case of Lemma C, carried out in Appendix of this present paper.

Let us investigate the boundary of concordance complement. First, note that the tubular neighbourhood $N(C)$ of the cylinder $C$ in $M \times I$ is\footnote{Fox and Milnor refer to a paper \cite{Noguchi} by Noguchi about surfaces in 4-space} of the form $C \times I^2$, where $I^2$ denotes a 2-cell, and $N(C)$ is homotopic to a circle $S^1$, this $S^1$ represents the class of the knot $k_0$ or $k_1$. Let $Q=\overline{M\times I - N(C)} $ denotes the closure of $M\times I - N(C)$ in $M\times I$. The intersection of $N(C)$ and $Q$ is $C\times \partial I^2$, that is homotopic to $S^1 \times S^1$, the first copy of $S^1$ represents $k_0$( or $k_1$) and the second copy represents a meridian of the corresponding knot $k_1$ (or $k_2$).
The boundary $\partial Q$ is the union of $C \times \partial I^2$ and the closure $Y_0,Y_1$ of $M \times \{0\} - N(C)$  and $M \times \{1\} - N(C)$ matched along tori $k_0 \times \partial I^2$ and $k_1 \times \partial I^2$ via zero surgeries $h_0, h_1$ correspondingly 
\[
\partial Q = Y_0 \cup_{h_0} C\times \partial I^2 \cup_{h_1} Y_1
\]

Note that, $Y_0$ and $Y_1$ are the closure of the knot complements in the homology 3-sphere $M$ for the knots $k_0,k_1$. That is
\begin{align*}
Y_0&=M-\text{Int} N(k_0),\\
Y_1&=M-\text{Int} N(k_1)
\end{align*}

where $N(k_0), N(k_1)$ are the tubular neighbourhood of $k_0,k_1$ in $M$.

Also the part $C \times \partial I^2$ can be imagine as the product of an annulus with a circle, so it looks like a solid torus with a torus tunnel inside along the core circle. So the boundary of $C \times \partial I^2$ has two component, that is two tori.

The zero surgeries $h_0,h_1$, gluing $Y_0,Y_1$ with $C \times \partial I^2$, are described as follows. The meridian of $\partial N(k_0)$ is glued along the longitude of one torus component $S^1 \times \{0\} \times \partial I^2$ of $C \times \partial I^2$, the meridian of $\partial N(k_1)$ is glued along the longitude of the other torus component $S^1 \times \{1\} \times \partial I^2$ of $C \times \partial I^2$.

Let $W_0$ be equals to $Y_0 \cup_{h_0} S^1 \times \{0\} \times \partial I^2$ and $W_1$ be equals to $Y_1 \cup_{h_1} C \times \partial I^2$. Then $\partial Q = W_0 \cup W_1$.

We will need the following Lemma. The proof of it is similar as of Kreinbihl (see \cite[Lemma 3.6]{Kreinbihl}).

{\bf Lemma 3}. Let $k$ be a knot in an oriented homology 3-sphere. Let N(k) be the tubular neighbourhood of $k$ in $M$. Let  $Y=M-\text{Int} N(k)$ be the closure complement space of $N(k)$ in $M$. Let $W=Y \cup_{h} S^1 \times S^1$ be the space, getting via a zero surgery $h$ on $Y$ along its boundary with a torus. Such that $h: S^1 \times S^1 \rightarrow \partial N(k)$ is the homeomorphism sending the longitude of $S^1 \times S^1$ to the meridian on $\partial N(k)$. Then, the homologies of $Y$ and $W$ are the same for all dimensions.

Proof.

We have $Y-\text{Int} (Y) = \partial N(k)$ and $W-\text{Int}(Y) = Z$, where $Z=\partial N(k) \cup_h S^1\times S^1$.  By Excision,

\[
H_*(Y,W)\cong H_*(\partial N(k),Z).
\]

In the long exact sequence of the pair $(\partial N(k),Z)$ the homomorphism $H_*(\partial N(k)) \rightarrow H_* Z$ is an isomorphism, sending fundamental class of $\partial N(k)$ to fundamental class of $Z$. Hence 

\[
H_*(\partial N(k),Z) \cong 0.
\]

Thus, $H_*(Y,W)\cong 0$. That means, $H_*(Y) \rightarrow H_*(W)$ in the long exact sequence of the pair $Y,W$ are isomorphisms. Q.E.D.

From the Lemma 3, we have the spaces $Y_0$ and $W_0$ have the same homology induced by inclusion $Y_0 \hookrightarrow W_0$. And it is also for $Y_1$ and $W_1$.

\section{Alexander polynomial, Alexander functions and Milnor torsion}

Let $k$ be a knot in an oriented homology 3-sphere $M$. There is a map $g: M-k \rightarrow S^1$ that is the extension of a map $\bar{g}: \partial (M-k) \rightarrow S^1$ (see \cite[p. 259]{AustinRolfsen}, \cite[p. 92]{Saveliev}). Then for $p \in S^1$ being a regular value, the closure $\overline{g^{-1}(p)}$ is an oriented surface $F$. And such that $k$ is the boundary $\partial F$ of $F$. The homotopy class $[g]$ in $[M-k,S^1] \cong H^1(M-k) \cong \mathbb{Z}$ \footnote{Homotpy construction of cohomology, see for example \cite[chapter 4]{Hatcher}} is a generator for $H^1(M-k)$. We have a pullback bundle

\begin{center}
\begin{tikzcd}
Y_{\infty} \arrow[r] \arrow[d] 
& \mathbb{R} \arrow[d] \\
M-k \arrow[r, "g"] & S^1 \\
\end{tikzcd}
\end{center}

where $Y_{\infty}$ is the infinite cyclic cover of $M-k$, defined by the kernel of Hurewicz map $\pi_1(M-k) \rightarrow \mathbb{Z} \cong H_1(M-k)$. The homology group $H_1(M-k)$ acts on $Y_{\infty}$ as group of deck transformations. Denoted by $t: Y_{\infty} \rightarrow Y_{\infty}$ a generator for $H_1(M-k)$. Then the first homology group of the cover $Y_{\infty}$ has a $\mathbb{Z}[t^{\pm}]$-module structure (see Milnor \cite{Milnor}). This module is called Alexander module. The order $\Delta_0(H_1(Y_{\infty}))$ of the module is the Alexander polynomial $\Delta_k(t)$ of the knot $k$, where $\Delta_k(t)$ is in the ring $\mathbb{Z}[t^{\pm}]$ defined up to a unit. 

The Alexander module can be described equivalently by homology with local coefficients\footnote{Detail on Homology with local coefficients can be found in the chapter 5 of a book by Davis and Kirk \cite{DavKir}. A direct application to define Alexander module in chapter 2 of Kreinbihl's PhD thesis \cite{Jamesthesis} at Wesleyan University, the thesis (available \url{https://digitalcollections.wesleyan.edu/object/ir-2227}) is a detailed version of the paper \cite{Kreinbihl}}, related to group ring $\mathbb{Z} [H_1(M-k)] \cong \mathbb{Z}[t^{\pm}]$. Then the Alexander module of the knot $k$ is homology $H_1(M-k;\mathbb{Z} [H_1(M-k)])$ with local coefficients in $\mathbb{Z} [H_1(M-k)]$. And the Alexander polynomial of the knot $k$ is the order of $H_1(M-k;\mathbb{Z} [H_1(M-k)])$
\[\Delta_k(t)=\Delta_0(H_1(M-k;\mathbb{Z} [H_1(M-k)]) = \Delta_0(H_1(Y_{\infty})) \in \mathbb{Z}[t^{\pm}] \]

We carry here definitions of Alexander function and its relation to Milnor torsion\footnote{See Turaev's survey on Reidemeister torsion in knot theory \cite{Turaev} for detail on Alexander function and Milnor torsion (section 1.1)} for knot complement $M-k$ in $M$. Let $\mathcal{K}$ is the field of quotients of $\mathbb{Z} [H_1(M-k)]$. Then the Alexander function $A(M-k)$ of the knot $M-k$ is defined to be is an element in $\mathcal{K}$, defined up to sign and multiplication by power $t$
\begin{align*}
A(M-k)&=\prod_{i=0}^m [\Delta_0(H_i(Y_{\infty}))]^{(-1)^{i+1}}\\
&= \prod_{i=0}^m [\Delta_0(H_i(M-k;\mathbb{Z} [H_1(M-k)]))]^{(-1)^{i+1}}.
\end{align*}

Now let us compute the Alexander function of a knot $k$ in a oriented homology 3-sphere $M$

{\bf Lemma 4.} Let $k$ be a knot in a oriented homology 3-sphere $M$. Then the Alexander function of $M-k$ is
\[
A(M-k)= \Delta_k(t)(t-1)^{-1}
\]
 Proof. 
 For every $i \geq 2$, the homology $H_i(M-k;\mathbb{Z}[t^{\pm}])$ is zero, so its orders are all 1.
 The zero dimensional chains $\mathbb{Z}[t^{\pm}]$-module $C_0(M-k,\mathbb{Z}[t^{\pm}]$ generated by a single point, that is a lift a a point in $M-k$, by the same argument\footnote{This argument is used severally times in the Milnor paper about Infinite cyclic covering, one can consult Kreinbihl's PhD thesis \cite[section 2.3.3, page 32]{Jamesthesis} for a detail description of the argument by using Fox calculus} using by Milnor (see \cite{Milnor}, Example 1, p. 119) we see that the order $\Delta_0(H_0(M-k;\mathbb{Z} [t^{\pm}]))$ is $(t-1)$. Obviously by definition the Alexander polynomial $\Delta_k(t)$ is the order $\Delta_0(H_1(M-k;\mathbb{Z} [t^{\pm}]))$. Q.E.D.

{\bf Theorem 5.} \cite[p. 126]{Turaev} Let $X$ is a CW complex then
$$\tau(X)~ \dot= ~A(X)$$
where $\tau(X)$ is the Milnor torsion of $X$, $A(X)$ is the Alexander function of $X$.

There is a property of torsion of pairs according to Kreinbihl (see \cite[Theorem 4.3.3]{Jamesthesis})

{\bf Theorem 6.} \cite{Jamesthesis} Let $\mathcal{M}$ be an $n$-manifold such that $\partial \mathcal{M} = M_1 \cup M_2$ and $M_1 \cap M_2 = \partial M_1 = \partial M_2$. Then,
\[
\tau(\mathcal{M},M_1) ~\dot= ~\overline{\tau(\mathcal{M},M_2)}^{(-1)^{n+1}}
\]
where the bar indicates the changing variable $t \rightarrow 1/t$.

\section{Proof of Theorem A}

We will prove Theorem A using properties of Milnor torsion and its relation to Alexander function as of Fox-Milnor \cite{FoxMilnor} and Kreinbihl \cite{Kreinbihl}. 

Let $k_0, k_1$ be concordant knots in an oriented homology 3-sphere $M$. Let $f$ be a concordance between $k_0$ and $k_1$ as described in section 3 of this paper.

We have seen by Lemma 3 that the homologies of the complement $Y_0$ of $k_0$ in $M$ and its corresponding zero surgery $W_0$ are the same. Thus we has the following equations between Alexander function of $Y_0$ and Milnor torsion of $W_0$ by Theorem 4
\[
A(Y_0)=A(W_0)~\dot= ~\tau(W_0);
\]
The similar is for the complement $Y_1$ of $k_1$
\[
A(Y_1)=A(W_1)~\dot= ~\tau(W_1);
\]
According to Milnor \cite[Lemma 4]{Milnor1}, we have the relation between torsion of the closure concordance complement $Q$ and $W_0$, recall that $W_0$ is a part of the boundary $\partial Q=W_0 \cup W_1$ of $Q$
\[
\tau(W_0)~\dot= ~\tau(Q,W_0)^{-1} \cdot \tau(Q)
\]
 and
 \[
\tau(Q)~\dot= ~\tau(Q,W_1) \cdot \tau(W_1)
\]
 Hence
 \[
 \tau(W_0)~\dot= ~\tau(Q,W_0)^{-1}\cdot \tau(Q,W_1) \cdot \tau(W_1)
 \]
 Note that, $Q,W_0,W_1$ satisfy the hypotheses of Theorem 6. So we have $ \tau(Q,W_1) ~\dot= ~\overline{\tau(Q,W_0)}^{(-1)}$. Hence
\[
\tau(W_0)~\dot= ~\tau(Q,W_0)^{-1}\cdot \overline{\tau(Q,W_0)}^{(-1)}  \cdot \tau(W_1)
\]

That means $A(Y_0) ~\dot= ~ \lambda(t) \lambda(1/t) A(Y_1)$, where $\lambda(t)$ denotes the rational function $\tau(Q,W_0)^{-1}$, that is an element of the field of fractions of  $\mathbb{Z}[t^{\pm}]$.

By Lemma 4 we have $A(Y_0)=\Delta_{k_0}(t)/(t-1)$ and $A(Y_1)=\Delta_{k_1}(t)/(t-1)$.
We conclude
\[\Delta_{k_0}(t) ~\dot{=} ~\lambda(t) \lambda(1/t) \Delta_{k_1}(t)\]
 Q.E.D.

{\bf Remark:} 
A particular case of concordance is sliceness of a knot, from Theorem A, if a knot in a oriented homology sphere is slice, by the same argument on the product $\lambda(t) \lambda(1/t)$ due to Milnor (see \cite{FoxMilnor}, Theorem 2, p. 263) then the Alexander polynomial of the knot is of the form $p(t)p(1/t)$, for $p(t)$ is a polynomial of integer coefficients in $t$. As I mentioned in the Note added in this version, an interesting proof of the Fox-Milnor Theorem has been established \cite{MarcPowell}.

\section{Proof of Theorem B}

The Theorem B is a direct generalisation from Theorem A. We consider the more general context, that is of the situation considered in Lemma C (see Appendix). This generalisation arise from a question due to Melikhov S.A. during a conversation with the author \cite{Melikhov}. 

Let $\mathcal{W}$ be homological cobordism between homological 3-spheres $M$ and $M'$. Let $k$ and $k'$ be knots in $M$ and $M'$ correspondingly. Let the concordance between $k$ and $k'$ be $g: S^1 \times I \rightarrow \mathcal{W}$. A tubular neighbourhood $N(C)$ of the 2-cell $C=g(S^1 \times I)$ in $\mathcal{W}$ is of the form $C\times I^2$ \footnote{cf. \cite{Noguchi}}, where $I^2$ denotes a 2-cell. The intersection of $N(C) \cap \partial \mathcal{W}$ of $N(C)$ with the boundary $\partial \mathcal{W}$ is the union tubular neighbourhoods $k\times I^2 \cup k' \times I^2$ of $k, k'$ in the corresponding homological spheres.

Consider the closure $Q$ of $\mathcal{W}- N(C)$ in the cobordism $\mathcal{W}$, the boundary $\partial Q$ of $Q$ is the union of $C \times \partial I^2$ and the closure of complements $Y,Y'$ of the knots $k,k'$ in $M,M'$ matched along the tori \footnote{The detail of this matching is similar to that of the surgery, described in section 3} $k \times \partial I^2$ and $k' \times \partial I^2$. Thus 
\[
\partial Q = Y \cup_{h} C\times \partial I^2 \cup_{h'} Y'
\]

where $h,h'$ are the corresponding zero surgery homeomorphisms, gluing $Y, Y'$ and $C\times \partial I^2$.

Denote $\mathcal{Y} = Y \cup_h S^1\times \{0\} \times \partial I^2$ and $\mathcal{Y'} = Y' \cup_{h'} C\times \partial I^2$. Rewrite $\partial Q = \mathcal{Y} \cup \mathcal{Y'}$.

By Lemma 3 and Theorem 4, we have
\[
A(Y)=A(\mathcal{Y})~\dot= ~\tau(\mathcal{Y});
\]
\[
A(Y')=A(\mathcal{Y'})~\dot= ~\tau(\mathcal{Y'});
\]
Analogy to the proof of Theorem A, we have
\[
\tau(\mathcal{Y})~\dot= ~\tau(Q,\mathcal{Y})^{-1} \cdot \tau(Q)
 \tau(\mathcal{Y})~\dot = ~\tau(Q,\mathcal{Y})^{-1}\cdot \tau(Q,\mathcal{Y'}) \cdot \tau(\mathcal{Y'}) \\ 
\]

As we see that $Q,\mathcal{Y},\mathcal{Y'}$ satisfy the hypotheses of Theorem 6.
Hence $\tau(Q, \mathcal{Y}) ~\dot= ~\overline{\tau(Q,\mathcal{Y'})}^{(-1)}$. Thus
\[
A(\mathcal{Y}) ~\dot= ~ \lambda(t) \lambda(1/t) A(\mathcal{Y'})
\]

or equivalently

\[\Delta_{k}(t) ~\dot{=} ~\lambda(t) \lambda(1/t) \Delta_{k'}(t)\]
for some polynomial $p(t)$. Q.E.D. 

\section*{Acknowledgements}

I would like to thank Prof. Kauffman and Prof. Melikhov for their engaging discussions. I am thankful to anonymous referee for his careful and supportive work. I am also grateful to everyone whose support and encouragement accompanied me along the way as I completed this paper.

The work was supported by the Ministry of Science and Higher Education of the Russian Federation (Agreement 075–02–2025–1728/2).

\section*{Appendix}

In this appendix we carry a more fact, stated in Lemma C of current appendix, due to S.A. Melikhov \cite{Melikhov}. The Lemma 2 in section 3 is a special case of the Lemma C.
Consider the following situation.
Let $\mathcal{W}$ be homological cobordism between homological 3-spheres $M$ and $M'$, i.e. the inclusions $M \hookrightarrow \mathcal{W}$, $M' \hookrightarrow \mathcal{W}$ induce isomorphisms $H_*(M) \rightarrow H_*(\mathcal{W})$, $H_*(M') \rightarrow H^*(\mathcal{W})$ on homology of $M,M'$ and $\mathcal{W}$. Then naturally the inclusions induce also isomorphisms $H^*(\mathcal{W}) \rightarrow H^*(M)$, $H^*(\mathcal{W}) \rightarrow H^*(M')$ on cohomology. The boundary $\partial \mathcal{W}$ is $M \cup M'$. Let $k$ be a knot in $M$ and $k'$ is a knot in $M'$. Suppose there is a concordance $g: S^1 \times I \rightarrow \mathcal{W}$ between $k = S^1 \times \{0\} \subset M$ and $k' = S^1 \times \{1\} \subset M'$ in $\mathcal{W}$. Denote by $C \subset \mathcal{W}$ the image of embedded cylinder $S^1 \times I$ under $g$. Let $N(C)$ be a regular neighbourhood of $C$ in $\mathcal{W}$. Denote by $Q$ the closure $\overline{\mathcal{W}-N(C)}$ of $\mathcal{W}-N(C)$.

{\bf Lemma C}. Concordant complement $\mathcal{W}-C$ has the homology of the circle $S^1$.

{\it Proof.} First, note that from the exact sequence of pair $(\mathcal{W}, M)$ we immediately have $H^i(\mathcal{W},M)=0$ for all dimension $i$.

Consider the cohomology exact sequence of triple $(\mathcal{W},\partial \mathcal{W}, M)$
\[
H^i(\mathcal{W},M) \rightarrow H^i(\partial \mathcal{W}, M) \rightarrow H^{i+1} (\mathcal{W},\partial \mathcal{W}) \rightarrow H^{i+1} (\mathcal{W},M)
\]

that is rewritten as
\[
0 \rightarrow H^i(\partial \mathcal{W}, M) \rightarrow H^{i+1} (\mathcal{W},\partial \mathcal{W}) \rightarrow 0
\]

Moreover, we have $H^i(\partial \mathcal{W}, M) = H^i(M')$ by Excision, then $H^{i+1}(\mathcal{W}, \partial \mathcal{W}) = H^i(M')$. By analogy we aslo have $H^{i+1}(\mathcal{W}, \partial \mathcal{W}) = H^i(M)$.

Now by Excision removing the complement of knots $k,k'$ in $M, M'$ correspondingly we have $H^i(N(C) \cup \partial \mathcal{W}, \partial \mathcal{W})=H^i(N(C), N(C) \cap \partial \mathcal{W})$.
Note that $N(C)$ is homotopic to $C$ and the intersection $N(C) \cap \partial \mathcal{W}$ is the disjoint union of the tubular neighbourhood of $k$ in $M$ with the tubular neighbourhood of $k'$ in $M'$. By homotopy invariance of cohomology we have $H^i(N(C), N(C) \cap \partial \mathcal{W}) = H^i(C, k \cup k')$. 

Consider the cohomology exact sequence of triple $(C, k \cup k', k')$:
\[
H^i(C,k')  \rightarrow H^i(k\cup k',k')  \rightarrow H^{i+1}(C, k \cup k')  \rightarrow H^{i+1}(C,k')
\]

Obviously $C$ and $k,k'$ are each homeomorphic to a circle, hence $H^i(C,k')=0$. We have $H^i(k\cup k',k')=H^(k)$ by Excision. Thus $H^{i+1}(C,k \cup k')=H^i(k)$.

From the above argument we conclude  $H^i(N(C) \cup \partial \mathcal{W}, \partial \mathcal{W})=H^{i-1}(k)$. 

Next we have the following equalities using first homotopy invariance of homology, then Poincar\'e-Lefschetz duality and Excision 
\[
H_{3-i}(\mathcal{W}-C)=H_{3-i}(Q)=H^{i+1} (Q,\partial Q) = H^{i+1}(\mathcal{W}, N(C) \cup \partial \mathcal{W}),
\]

where $\partial Q$ is the boundary of $Q$ and $\partial \mathcal{W} = M \cup M'$ is the boundary of $\mathcal{W}$. 

The cohomology groups $ H^{i+1}(\mathcal{W}, N(C) \cup \partial \mathcal{W})$ are included in the cohomology exact sequence of triple $(\mathcal{W}, N(C) \cup \partial \mathcal{W}, \partial \mathcal{W})$.
\[
H^i (\mathcal{W},\partial \mathcal{W})  \rightarrow H^i(N(C) \cup \partial \mathcal{W}, \partial \mathcal{W})  \rightarrow H^{i+1}(\mathcal{W}, N(C) \cup \partial \mathcal{W})  \rightarrow H^{i+1} (\mathcal{W}, \partial \mathcal{W})
\]

Replace $H^i (\mathcal{W},\partial \mathcal{W})$ by $H^{i-1}(M)$, $H^{i+1}(\mathcal{W}, N(C) \cup \partial \mathcal{W})$ by $H_{3-i}(\mathcal{W}-C)$ and $H^i(N(C) \cup \partial \mathcal{W}, \partial \mathcal{W})$ by $H^{i-1}(k)$ we have
\[
H^{i-1}(M)  \rightarrow H^{i-1}(k)  \rightarrow H_{3-i}(\mathcal{W}-C)  \rightarrow H^{i} (M)
\]

Clearly $H_4(\mathcal{W}-C)=0$, for $i=2$ we get $H_1(\mathcal{W}-C) = \mathbb{Z}$, for $i=3$ we get $H_0(\mathcal{W}-C)=\mathbb{Z}$. Also we have $H^0(M) \rightarrow H^0(k)$ is an isomorphism, hence $H_2(\mathcal{W}-C)=H_3(\mathcal{W}-C)=0$. Q.E.D.

Another proof for Lemma 1.

{\bf Lemma 1}. Let $k$ be a knot in an oriented homology 3-sphere $M$. Then the knot complement $M-k$ in $M$ has the homology of a circle, that is $H_*(M-k)=H_*(S^1)$.

{\it Proof.}
Let $\nu(k)$ be a tubular neighbourhood of the knot $k$ in the homology 3-sphere $M$. Denote by $X_k$ the closure of $M-\nu(k)$. By homotopy invariance of homology we have $H_{3-i}(M-k) = H_{3-i}(X_k)$, where $i \in \{ 0,1,2,3\}$. The Poincar\'e-Lefschetz duality and Excision give

\[
H_{3-i}(M-k) = H_{3-i}(X_k) = H^i(X_k, \partial X_k) = H^i(M,\nu(k))
\]

The groups $H^i(M,\nu(k)$ are included in the cohomology exact sequence of the pair $(M, \nu(k))$.

\[
H^{i-1}(M) \rightarrow H^{i-1}(\nu(k)) \rightarrow H^i (M,\nu(k)) \rightarrow H^i(M)
\]

or equivalently 

\[
H^{i-1}(M) \rightarrow H^{i-1}(k) \rightarrow H_{3-1} (M-k) \rightarrow (M)
\]

From the above exact sequence, let $i=2$ we get $H_1(M-k) = \mathbb{Z}$ and let $i=3$ we get $H_0(M-k) = \mathbb{Z}$. Clearly $H^0(M) \rightarrow H^0(k)$ is an isomorphism, then $H_2(M-k)=H_3(M-k)=0$. Q.E.D.

\end{document}